\documentclass[12pt]{article}
\usepackage{amssymb}

\usepackage{amsmath,latexsym,amsfonts}
\usepackage{graphicx}
\usepackage{srcltx}
\usepackage{amscd}

\newtheorem{lemma}{Lemma}[section]
\newtheorem{coro}[lemma]{Corollary}
\newtheorem{prop}[lemma]{Proposition}
\newtheorem{thm}[lemma]{Theorem}

\newtheorem{rem}[lemma]{Remark}
\makeatletter\@addtoreset{equation}{section}
\renewcommand\theequation{\thesection.\@arabic\c@equation}

\begin{document}
\begin{center}
{\LARGE Vanishing theorems  of the basic harmonic forms on a complete foliated Riemannian manifold}

 \renewcommand{\thefootnote}{}
\renewcommand{\thefootnote}{\arabic{footnote}}
\setcounter{footnote}{0}

\vspace{0.5 cm} {\large Seoung Dal Jung and Huili Liu}
\end{center}
\vspace{0.5cm}
\noindent {\bf Abstract.}  It is well-known that on a compact foliated Riemannian  manifold $(M,\mathcal F)$ with some transversal curvature conditions,  there are no nontrivial basic harmonic $r\ (0<r<q={\rm codim}\mathcal F)$-forms ( M. Min-Oo et al., J. Reine Angew. Math. 415 (1991) [\ref{MO}]).   In this paper, we extend the above facts to a complete foliated Riemannian
manifold.  
\section{Introduction}
Let $(M,g,\mathcal F)$ be a  foliated Riemannian manifold with a 
 foliation $\mathcal F$ and a bundle-like metric $g$ with respect to $\mathcal F$. A foliated Riemannian manifold is a Riemannian manifold with a Riemannian foliation, i.e., a foliation on a smooth manifold such that the normal bundle  is endowed with a metric whose Lie derivative is zero along leaf directions (see [\ref{RE}]). A Riemannian metric on $M$ is  bundle-like if the leaves of the foliation $\mathcal F$ are locally equidistant, that is, the metric $g$ on $M$ induces a holonomy invariant transverse metric on the normal bundle $Q=TM/T\mathcal F$, where $T\mathcal F$ is the tangent bundle of $\mathcal F$. Every Riemannian foliation admits bundle-like metrics. Many researchers have studided basic forms and the basic Laplacian on foliated Riemannian manifolds. Basic forms are locally forms on the space of leaves; that is, forms  $\phi$ satisfying $i(X)\phi=i(X)d\phi=0$ for all $X\in T\mathcal F$. Basic forms are preserved by the exterior derivative and are used to define basic de-Rham cohomology groups $H_B^* (\mathcal F)$. The basic Laplacian $\Delta_B$ for a given bundle-like metric is a version of the Laplace operator that preserves the basic forms. It is well-known [\ref{KT2},\ref{TO}] that on a closed oriented manifold $M$ with a transversally oriented Riemannian foliation $\mathcal F$, $H_B^r (\mathcal F)\cong \mathcal H_B^r (\mathcal F)$, where $\mathcal H_B^r(\mathcal F)={\rm ker}\Delta_B$ is finite dimensional. And so $\chi_B(\mathcal F)=\sum_{r=0}^q (-1)^r \dim \mathcal H_B^r (\mathcal F)$, where $\chi_B(\mathcal F)$ is the basic Euler characteristic [\ref{BPR}]. In 1991, M. Min-Oo et al. [\ref{MO}] proved that on a closed foliated Riemannian manifold $M$, if the transversal  curvature operator of $\mathcal F$ is positive definite, then $H_B^r (\mathcal F)=0 \ (0<r<q)$, that is, any basic harmonic $r$-form is trivial.

  In this paper, we study the basic $r$-forms on a complete foliated Riemannian manifold.
  
\bigskip
\noindent{\bf Main Theorem.} {\it Let $(M,g,\mathcal F)$ be a complete foliated Riemannian manifold and all leaves be compact. Assume that the mean curvature form is bounded and coclosed.

\noindent $(1)$ If the transversal Ricci cuvature of $\mathcal F$ is positive-definite, then  any $L^2$-basic harmonic $1$-forms $\phi$ with $\phi\in \mathcal S_B$ are trivial. 

\noindent$(2)$ If the curvature endomorphism of $\mathcal F$ is positive-definite, then any $L^2$-basic harmonic $r$-forms $\phi$ with $\phi\in \mathcal S_B$ are trivial.

\bigskip
\noindent Here $\mathcal S_B$ is the Sobolev space of basic forms whose derivative belong to $L^2\Omega_B^*(\mathcal F)$.
 }

\bigskip
\noindent Note that in 1980, H. Kitahara [\ref{KI1}] proved that under the same condition of the transversal Ricci curvature, there are no nontrivial basic $\Delta_T$-harmonic 1-forms with finite global norms. Here $\Delta_T$ is a different operator to the basic Laplacian $\Delta_B$. If $\mathcal F$ is minimal, then $\Delta_T = \Delta_B$.

\section{Preliminaries}
Let $(M,g,\mathcal F)$ be a $(p+q)$-dimensional complete foliated Riemannian
manifold with a foliation $\mathcal F$ of codimension $q$ and a
bundle-like metric $g$ with respect to $\mathcal F$. 
Let $TM$ be the tangent bundle of $M$, $T\mathcal F$ its integrable
subbundle given by $\mathcal F$, and  $Q=TM/T\mathcal F$ the corresponding
normal bundle of $\mathcal F$.
 Then we have an exact sequence of vector bundles
\begin{align}\label{eq1-1}
 0 \longrightarrow T\mathcal F \longrightarrow
TM_{\buildrel \longleftarrow \over \sigma }^{\buildrel \pi \over
\longrightarrow} Q \longrightarrow 0,
\end{align}
where $\pi:TM\to Q$ is a projection and $\sigma:Q\to T\mathcal F^\perp$ is a bundle map satisfying
$\pi\circ\sigma=id$.  
    Let $g_Q$ be the holonomy invariant metric
on $Q$ induced by $g$, i.e., $\theta(X)g_Q=0$ for any  vector field $X\in T\mathcal F$, where
$\theta(X)$ is the transverse Lie derivative [\ref{KT1}]. Let $R^Q$ and ${\rm Ric}^Q$ be the transversal curvature tensor and transversal Ricci operator of $\mathcal F$ with respect to the transversal Levi-Civita connection $\nabla^Q\equiv \nabla$ in $Q$ [\ref{TO}], respectively.
A differential form $\phi\in \Omega^r(M)$ is {\it basic} if $
i(X)\phi=0$ and $i(X)d\phi=0$ for all $X\in T\mathcal F$. 
In a distinguished chart $(x_1,\cdots,x_p; y_1, \cdots,y_q)$ of $\mathcal F$, a basic $r$-form $\phi$ is expressed by
\begin{align*}
\phi=\sum_{a_1 <\cdots<a_r} \phi_{a_1\cdots a_r}dy_{a_1}\wedge\cdots\wedge dy_{a_r},
\end{align*}
  where  the functions $\phi_{a_1\cdots a_r}$ are independent of $x$.
Let $\Omega_B^r(\mathcal F)$ be the set of all basic $r$-forms on
$M$. Then $\Omega^r(M)=\Omega_B^r(\mathcal F)\oplus \Omega_B^r(\mathcal F)^\perp$ [\ref{LO}].
  Now, we recall the star operator $\bar *:\Omega_B^r (\mathcal F)\to \Omega_B^{q-r}(\mathcal F)$ given by [\ref{KT2},\ref{PR}]
\begin{align}
\bar *\phi= (-1)^{p(q-r)}*(\phi\wedge\chi_{\mathcal F}),\quad\forall \phi\in\Omega_B^r(\mathcal F),
\end{align}
where $\chi_{\mathcal F}$ is the characteristic form of $\mathcal F$ and $*$ is the Hodge star operator associated to $g$.
For any basic forms $\phi,\psi \in \Omega_B^r(\mathcal F)$, it is well-known [\ref{PR}] that $\phi\wedge\bar *\psi = \psi\wedge\bar *\phi$ and $\bar *^2\phi = (-1)^{r(q-r)}\phi$.  
 The operator $d_B$ is the restriction of $d$ to the basic forms, i.e., $d_B= d|_{\Omega_B^*(\mathcal F)}$.  Let  $d_T = d_B-\kappa_B\wedge $ and $\delta_T = (-1)^{q(r+1)+1}\bar * d_B\bar *$, where  $\kappa_B$ is the basic part of the mean curvature form $\kappa$ of $\mathcal F$ [\ref{LO}].  Note that $\kappa_B$ is closed, i.e., $d\kappa_B=0$ [\ref{PJ},\ref{TO}].
The operator $\delta_B:\Omega_B^r(\mathcal F)\to \Omega_B^{r-1}(\mathcal F)$  is defined by
\begin{align}
\delta_B\phi = (-1)^{q(r+1)+1}\bar * d_T\bar *\phi = \delta_T\phi + i(\kappa_B^\sharp)\phi,
\end{align}
where  $(\cdot)^\sharp$ is the $g_Q$-dual vector field of $(\cdot)$.    Generally,  $\delta_B$ is not a restriction of $\delta$ on $\Omega_B^r (\mathcal F)$, i.e., $\delta_B \ne \delta|_{\Omega_B^r (\mathcal F)}$, where $\delta$ is the formal adjoint of $d$. But $\delta_B\omega=\delta\phi$ for any basic $1$-form $\phi$. 
  Let $\Delta_B =d_B\delta_B +\delta_B d_B$ be a basic Laplacian. Then  $\Delta^M|_{\Omega_B^0(\mathcal F) } =\Delta_B$ [\ref{KT2}], where $\Delta^M$ is the Laplacian on $M$. 
 Let
 $\{E_a\}(a=1,\cdots,q)$ be a local orthonormal basic  frame of $Q$ and $\theta^a$ a $g_Q$-dual $1$-form to $E_a$.  We define  $\nabla_{\rm tr}^*\nabla_{\rm tr}:\Omega_B^r(\mathcal F)\to \Omega_B^r(\mathcal F)$ by
\begin{align}\label{eq1-12}
\nabla_{\rm tr}^*\nabla_{\rm tr} =-\sum_a \nabla^2_{E_a,E_a}
+\nabla_{\kappa_B^\sharp},
\end{align}
where $\nabla^2_{X,Y}=\nabla_X\nabla_Y -\nabla_{\nabla^M_XY}$ for
any $X,Y\in TM$ and $\nabla^M$ is the Levi-Civita connection with respect to $g$.     
Then the  generalized Weitzenb\"ock type formula on $\Omega_B^r(\mathcal F)$ is given by [\ref{Jung}]
\begin{align}\label{2-7}
\Delta_B \phi = \nabla_{\rm tr}^*\nabla_{\rm tr} \phi
+ F(\phi) + A_{\kappa_{B}^\sharp} \phi
\end{align}
for any $\phi\in\Omega_B^r(\mathcal F)$,
where  $F=\sum_{a,b=1}^{q}\theta^a\wedge i(E_b)  R^Q(E_b,E_a)$ and 
\begin{align}\label{eq1-13}
A_Y\phi =\theta(Y)\phi-\nabla_Y\phi.
\end{align}
In particular, for a 1-form $\phi$, $F(\phi)^\sharp = {\rm Ric}^Q (\phi^\sharp)$ and $A_Y s = -\nabla_{\sigma(s)}\pi(Y)$. 
Let $\Omega_{B,o}^*(\mathcal F)$ be the space of basic forms with compact supports.

Let $\nu$ be the transversal volume form, i.e., $*\nu =\chi_{\mathcal F}$. The pointwise inner product   $\langle\cdot,\cdot\rangle$ on $\Omega_B^r(\mathcal F)$ is given by
\begin{align}
\langle\phi,\psi\rangle \nu = \phi\wedge\bar *\psi
\end{align}
for any basic forms $\phi,\psi \in \Omega_B^r(\mathcal F)$. And the global inner product   $\ll \cdot,\cdot\gg$  on $\Omega_B^*(\mathcal F)$ is defined by
\begin{align}
\ll \phi,\psi\gg =\int_M \langle\phi,\psi\rangle\mu_M 
\end{align}
for any $\phi,\psi\in\Omega_B^r(\mathcal F)$, one of which has compact support, where  $\mu_M=\nu\wedge\chi_{\mathcal F}$ is the volume form with respect to $g$.  It is well-known [\ref{Jung}] that $\ll \nabla_{\rm tr}^*\nabla_{\rm tr}\phi,\psi\gg =\ll\nabla_{\rm tr}\phi,\nabla_{\rm tr}\psi\gg$ for any $\phi,\psi\in\Omega_{B,0}^r(\mathcal F)$ and
\begin{align}
\ll d_B\phi,\psi\gg &=\ll\phi,\delta_B\psi\gg
\end{align}  
for any $\phi\in\Omega_{B,o}^r(\mathcal F), \psi\in\Omega_{B,o}^{r+1}(\mathcal F)$. 
The basic form $\phi$ is said to be {\it $L^2$-basic form} if $\phi$ has finite global norm, i.e., $\Vert\phi\Vert^2<\infty$.  Let $\mathcal H_{B,2}^r(\mathcal F)$ be the space of $L^2$-basic harmonic forms, i.e., 
\begin{align}
\mathcal H_{B,2}^r(\mathcal F)=\{\phi\in L^2\Omega_B^r(\mathcal F)\ |\ d_B\phi =\delta_B\phi =0\}.
\end{align}
 Generally, the space $\mathcal H_{B,2}^r(\mathcal F)$ can have infinite dimension. And if the dimension of $\mathcal H_{B,2}^r(\mathcal F)$ is finite, then it depends on the bundle-like metric. Trivially, if $M$ is compact, then $\mathcal H_{B,2}^r (\mathcal F)\cong   H_B^r(\mathcal F)$. So we study the vanishing properties of the $L^2$-basic harmonic spaces on a complete foliated Riemannian manifold. 
\begin{rem} {\rm (1) The operator $\delta_T$ is the formal adjoint of $d_B$ with respect to the global norm $(\cdot,\cdot)$, which is given  by
\begin{align}
(\phi,\psi)=\int_M \langle\phi,\psi\rangle \nu\wedge dx_1\wedge\cdots\wedge dx_p.
\end{align}
Let $\Delta_T = d_B\delta_T +\delta_T d_B$ be a Laplacian. If the foliation is  minimal, then $\delta_T=\delta_B$. So $\Delta_B =\Delta_T$.

 (2)  
 In 1980, H. Kitahara [\ref{KI1}] proved that if the transversal Ricci curvature is nonnegative and positive at some point, then there are no nontrivial $L^2$-basic $\Delta_T$-harmonic $1$-forms.
 }
 \end{rem}

\section{Vanishing theorem}
Let $(M,g,\mathcal F)$ be a complete foliated Riemannian manifold with a foliation $\mathcal F$ of codimension $q$ and a bundle-like metric $g$ with respect to $\mathcal F$. Assume that all leaves of $\mathcal F$ are compact.
Now, we consider a smooth function $\mu$ on $\mathbb R$ satisfying
\begin{align*}
(i) \ 0\leq \mu(t)\leq 1\ {\rm on}\ \mathbb R,\quad
(ii)\ \mu(t)=1\ \ {\rm for}\ t\leq 1,\quad
(iii)\ \mu(t)=0\ \ {\rm for}\ t\geq 2.
\end{align*}
Let $x_0$ be a point in $M$. For each point $y \in M$, we
denote by $\rho(y)$ the distance between leaves through $x_0$ and $y$. For any real number $l>0$, we define a Lipschitz continuous function $\omega_l$  on $M$ by  
\begin{align*}
\omega_l(y)=\mu(\rho(y)/l).
\end{align*} 
Trivially, $\omega_l$ is a basic function. 
   Let $B(l) = \{ y \in M | \rho(y) \leq l \}$ for $l>0$. Then $\omega_l$ satisfies
the following properties:
\begin{center}
 $\left.
 \begin{array}{ll}
0 \leq \omega_l(y) \leq 1 & \text{for any} \ y \in M\\
\text{supp}\ \omega_l \subset B(2l) & \\
\omega_l(y)=1 & \text{for any} \  y \in B(l)\\
\lim_{l \rightarrow \infty} \omega_l = 1 & \\
|d_B\omega_l | \leq \frac{C}{l} & \text{almost everywhere on}\ M, 
\end{array} \right.$
\end{center}
where $C$ is a positive constant independent of $l$ [\ref{YO1}]. Hence $\omega_l\psi$ has compact support for any basic form $\psi\in\Omega_B^*(\mathcal F)$ and $\omega_l\psi \to \psi$ (strongly) when $l\to \infty$.
\begin{lemma} $[\ref{KI}]$ For any $\phi\in\Omega_B^r (\mathcal F)$,  there exists a number $A$ depending only on $\mu$, such that
\begin{align*}
&\Vert d_B\omega_l \wedge\phi \Vert_{B(2l)}^2\leq {qA^2 \over l^2} \Vert\phi \Vert_{B(2l)}^2,\\
&\Vert d_B\omega_l \wedge \bar * \phi\Vert_{B(2l)}^2 \leq {qA^2\over l^2}\Vert \phi\Vert_{B(2l)}^2, \\
&\Vert d_B\omega_l \otimes\phi \Vert_{B(2l)}^2\leq {qA^2 \over l^2} \Vert\phi \Vert_{B(2l)}^2,
\end{align*}
where $\Vert \phi \Vert_{B(2l)}^2 =\ll\phi,\phi\gg_{B(2l)}=\int_{B(2l)} \langle\phi,\phi\rangle\mu_M$.
\end{lemma}

\begin{prop}  For any $L^2$-basic form $\psi$,   if $\Delta_B\psi=0$,  then $d_B\psi=0$ and $\delta_B\psi=0$.
\end{prop}
{\bf Proof.} Let $\psi$ be a $L^2$-basic form. Then we have
\begin{align}\label{3-1}
 \ll \Delta_B\psi,\omega_l^2 \psi\gg_{B(2l)} =\ll d_B\psi,d_B(\omega_l^2\psi)\gg_{B(2l)} +\ll\delta_B\psi,\delta_B(\omega_l^2\psi)\gg_{B(2l)}.
\end{align} 
 By a direct calculation, we have
\begin{align}
d_B(\omega_l^2\psi) &=   \omega_l^2 d_B\psi +2\omega_l d_B\omega_l \wedge\psi,\\
\delta_B (\omega_l^2\psi) &= \omega_l^2 \delta_B \psi + (-1)^{q(r+1)+1}\bar * (2\omega_l d_B \omega_l \wedge\bar * \psi).
\end{align}
From (3.1), (3.2) and (3.3), if $\Delta_B\psi=0$, then      
\begin{align*}
&\Vert\omega_l d_B\psi\Vert^2_{B(2l)} + \Vert\omega_l \delta_B\psi\Vert^2_{B(2l)} \\
&=-2\ll \omega_ld_B\psi,d_B\omega_l\wedge\psi\gg_{B(2l)}  
+ 2(-1)^{q(r+1)}\ll\omega_l\delta_B\psi,\bar *(d_B\omega_l\wedge \bar *\psi)\gg_{B(2l)}.
\end{align*}
Hence  by the  the Schwartz's inequality and Lemma 3.1, we have
\begin{align*}
&\Vert\omega_l d_B\psi\Vert^2_{B(2l)} + \Vert\omega_l \delta_B\psi\Vert^2_{B(2l)}\\
&\leq {\epsilon_1}\Vert \omega_l d_B\psi\Vert^2_{B(2l)} +{\epsilon_2}\Vert\omega_l \delta_B \psi \Vert^2_{B(2l)}
+{B_1\over l} \Vert\psi\Vert^2_{B(2l)}
\end{align*}
for some positive real numbers $\epsilon_1,\epsilon_2$ and $B_1$.  Therefore, we have
\begin{align*}
\Vert\omega_l d_B\psi\Vert^2_{B(2l)} + \Vert\omega_l \delta_B\psi\Vert^2_{B(2l)} \leq {B_2\over l} \Vert \psi\Vert_{B(2l)}^2
\end{align*}
for some positive real number $B_2$. Since $\psi$ is the $L^2$-basic form, letting $l\to \infty$, $d_B\psi=\delta_B\psi=0$. $\Box$
\begin{rem} {\rm  In 1979, H. Kitahara [\ref{KI}] proved the corresponding result with the Laplacian  $\Delta_T$. Namely, on a complete foliated manifold, if $\Delta_T\phi=0$, then $d_B\phi =\delta_T\phi=0$.   }
\end{rem}  
  Now we prove  the vanishing theorem of the $L^2$-basic harmonic form on a complete foliated Riemannian manifold. First of all, we prepare some lemmas.
  \begin{lemma} Let $(M,g,\mathcal F)$ be a complete foliated Riemannian manifold whose leaves are compact.
  Suppose that $\kappa_B$ is bounded and coclosed. Then for any $L^2$-basic harmonic form $\phi$,
  \begin{align}
 \limsup_{l\to\infty} \ll A_{\kappa_B^\sharp}\phi,\omega_l^2\phi\gg_{B(2l)}=0.
  \end{align}
  \end{lemma}
  {\bf Proof.} Let $\phi$ be a  $L^2$-basic harmonic form. Since $\theta(X)\phi=d_B i(X)\phi$,  from (2.6)  we have
 \begin{align}
 \ll A_{\kappa_B^\sharp}\phi,\omega_l^2\phi\gg_{B(2l)} =\ll d_B i(\kappa_B^\sharp)\phi,\omega_l^2\phi\gg_{B(2l)} -\ll\nabla_{\kappa_B^\sharp}\phi,\omega_l^2\phi\gg_{B(2l)}
 \end{align}
 Since $\delta_B\phi=0$, from (3.3) and Lemma 3.1, we have
 \begin{align*}
| \ll d_B i(\kappa_B^\sharp)\phi,\omega_l^2\phi\gg_{B(2l)}|&= 2|\ll \omega_l i(\kappa_B^\sharp)\phi, \bar *(d_B\omega_l\wedge \bar *\phi)\gg_{B(2l)}|\\
 &\leq \epsilon_3 \Vert\omega_l i(\kappa_B^\sharp)\phi\Vert^2_{B(2l)}+{1\over\epsilon_3} \Vert d_B\omega_l\wedge \bar *\phi\Vert^2_{B(2l)}\\
 &\leq \epsilon_3 \Vert\omega_l i(\kappa_B^\sharp)\phi\Vert^2_{B(2l)}+ {B_3\over l^2} \Vert\phi\Vert^2_{B(2l)}
 \end{align*} 
 for some positive real numbers $\epsilon_3$ and $B_3$.
 By using  $|i(\kappa_B^\sharp)\phi|^2 + |\kappa_B\wedge\phi|^2 =|\kappa_B|^2 |\phi|^2$, we have
  \begin{align}
| \ll d_B i(\kappa_B^\sharp)\phi,\omega_l^2\phi\gg_{B(2l)}|\leq \epsilon_3\ {\rm max}(|\kappa_B|^2) \Vert\omega_l\phi\Vert^2_{B(2l)} +{B_3\over l^2} \Vert\phi\Vert^2_{B(2l)}
\end{align}
On the other hand, since $\delta_B\kappa_B=0$, by a direct calculation, we have
\begin{align*}
\ll \nabla_{\kappa_B^\sharp}\phi,\omega_l^2\phi\gg_{B(2l)} &=\frac12 \ll d_B(|\omega_l\phi|^2),\kappa_B\gg_{B(2l)} - \ll \omega_l\phi,\kappa_B^\sharp(\omega_l)\phi\gg_{B(2l)}\\
&=-\ll \omega_l\phi,\kappa_B^\sharp(\omega_l)\phi\gg_{B(2l)}.
\end{align*}
Hence by the Schwartz inequality, we have
  \begin{align}
 | \ll \nabla_{\kappa_B^\sharp}\phi,\omega_l^2\phi\gg_{B(2l)}| &=  |\ll \omega_l\phi,\kappa_B^\sharp(\omega_l)\phi\gg_{B(2l)}|\notag\\
 &\leq \epsilon_4\Vert\omega_l\phi\Vert^2_{B(2l)} +{B_4\over l^2}{\rm max}(|\kappa_B|^2) \Vert\phi\Vert^2_{B(2l)}
  \end{align}
  for a positive real numbers $\epsilon_4$ and $B_4$.
From (3.6) and (3.7),  by letting $l\to \infty$, we have 
\begin{align*}
&\limsup_{l\to \infty} | \ll d_B i(\kappa_B^\sharp)\phi,\omega_l^2\phi\gg_{B(2l)}|\leq \epsilon_3\ {\rm max}(|\kappa_B|^2) \Vert\phi\Vert^2,\\
&\limsup_{l\to\infty}| \ll \nabla_{\kappa_B^\sharp}\phi,\omega_l^2\phi\gg_{B(2l)}| \leq \epsilon_4\Vert\phi\Vert^2. 
\end{align*}
Since $\epsilon_3$ and $\epsilon_4$ are arbitrary positive numbers, we have
\begin{align}
&\limsup_{l\to \infty} | \ll d_B i(\kappa_B^\sharp)\phi,\omega_l^2\phi\gg_{B(2l)}|=0, \\
&\limsup_{l\to\infty}| \ll \nabla_{\kappa_B^\sharp}\phi,\omega_l^2\phi\gg_{B(2l)}| =0. 
\end{align}
Hence from (3.5), (3.8) and (3.9), the proof is completed. $\Box$
  
  \begin{thm} Let $(M,g,\mathcal F)$ be as in Lemma 3.4.
  Suppose that $\kappa_B$ is bounded and coclosed.  If the curvature endomorphism $F$ of $\mathcal F$ is positive-definite, then any $L^2$- basic harmonic r-forms $\phi$ with $\phi\in \mathcal S_B$ are trivial, i.e., $\mathcal H_{B,2}^r(\mathcal F)=\{0\}$.
  \end{thm}
  {\bf Proof.} Let $\phi$ be a $L^2$-basic harmonic $r$-form.  From (2.5) and Proposition 3.2, we have
  \begin{align}
  \langle\nabla_{\rm tr}^*\nabla_{\rm tr}\phi,\omega_l^2 \phi\rangle +\langle F(\phi),\omega_l^2\phi\rangle + \langle A_{\kappa_B^\sharp}\phi,\omega_l^2\phi\rangle =0.
  \end{align}
 On the other hand, a direct calculation gives
 \begin{align}
 \ll \nabla_{\rm tr}^*\nabla_{\rm tr}\phi,\omega_l^2\phi \gg_{B(2l)}  &=\ll\nabla_{\rm tr}\phi,2\omega_l d_B\omega_l\otimes \phi \gg_{B(2l)} + \Vert\omega_l \nabla_{\rm tr}\phi\Vert^2_{B(2l)}.
 \end{align}
 From Lemma 3.1, we have
 \begin{align*}
| \ll\nabla_{\rm tr}\phi,2\omega_l d_B\omega_l\otimes \phi \gg_{B(2l)}|\leq {\epsilon_5}\Vert\omega_l\nabla_{\rm tr}\phi\Vert^2_{B(2l)} + {B_5\over l^2}\Vert\phi\Vert^2_{B(2l)}
 \end{align*}
 for some positive constants $\epsilon_5$ and $ B_5$. Hence by letting $l\to\infty$, we have
 \begin{align*}
\limsup_{l\to\infty} \ll\nabla_{\rm tr}\phi,2\omega_l d_B\omega_l\otimes \phi \gg_{B(2l)}\leq {\epsilon_5}\Vert\nabla_{\rm tr}\phi\Vert^2.
 \end{align*}
Since $\epsilon_5$ is arbitrary and $\phi\in \mathcal S_B$ ( i.e., $\Vert\nabla_{\rm tr}\phi\Vert^2<\infty$), we have
  \begin{align}
\limsup_{l\to\infty} \ll\nabla_{\rm tr}\phi,2\omega_l d_B\omega_l\otimes \phi \gg_{B(2l)}=0.
 \end{align}
Hence from  (3.11), (3.12) and Lemma 3.4, we have
 \begin{align}
 \Vert \nabla_{\rm tr}\phi\Vert^2  + \limsup_{l\to\infty} \ll F(\phi),\omega_l^2\phi\gg_{B(2l)} =0,
 \end{align}
 which complete the proof. $\Box$
  
  Since $F(\phi^\sharp)={\rm Ric}^Q(\phi^\sharp)$ for any basic 1-form $\phi$, we have the following corollary.
  \begin{coro} Let $(M,g,\mathcal F)$ be as in Lemma 3.4. Suppose that $\kappa_B$ is bounded and coclosed. If the transversal Ricci curvature ${\rm Ric}^Q$ is positive-definite, then any $L^2$-basic harmonic $1$-forms $\phi$ with $\phi\in \mathcal S_B$ are trivial, $\mathcal H_{B,2}^1(\mathcal F)=\{0\}$.
  \end{coro}


\bigskip
\noindent{\bf Acknowledgements.}  The first author was supported by  the National Research Foundation of Korea(NRF) grant funded
       by the Korea government (MSIP) (NRF-2015R1A2A2A01003491) and the second author was supported by NSFC (No. 11371080).

\noindent Department of Mathematics and Research Institute for Basic Sciences, Jeju National University, Jeju 690-756, Korea

\noindent {\it E-mail address} : sdjung@jejunu.ac.kr

\bigskip
\noindent Department of Mathematics, Northeastern University, Shenyang 110004, P. R. China

\noindent {\it E-mail address} : liuhl@mail.neu.edu.cn


\begin{thebibliography}{[9]}

\bibitem{Lop}\label{LO} J. A. Alvarez {L\'opez},
\emph{The basic component of the mean curvature of Riemannian
foliations}, Ann. Global Anal. Geom. 10 (1992), 179-194.



\bibitem{BPR}\label{BPR} V. Blefi, E. Park and K. Richardson, \emph{A Hopf index theorem for foliations}, Diff. Geom. Appl. 18 (2003), 319-341.

\bibitem{He}\label{HE} J. J. Hebda, \emph{Curvature and focal points in Riemannian foliations}, Indiana Univ. Math. J. 35 (1986), 321-331.




\bibitem{Jung}\label{Jung} S. D. Jung,
\emph{The first eigenvalue of the transversal Dirac operator}, J.
Geom. Phys. 39 (2001), 253-264.



\bibitem{Kamber2}\label{KT1} F. W. Kamber and Ph.
Tondeur, \emph{Infinitesimal automorphisms and second variation of
the energy for harmonic foliations}, T\^ohoku Math. J. 34 (1982),
525-538.

\bibitem{jj}\label{KT2} F. W. Kamber and Ph. Tondeur, {\it De Rham-Hodge theory for Riemannian foliations}, Math. Ann. 277 (1987), 415-431.

\bibitem{k}\label{KI} H. Kitahara, {\it Remarks on square-integrable basic cohomology spaces on a foliated Riemannian manifold}, Kodai Math. J. 2 (1979), 187-193.

\bibitem{kk}\label{KI1} H. Kitahara, {\it Nonexistence of nontrivial $\square^{\prime\prime}$-harmonic 1-forms on a complete foliated Riemannian manifold}, Trans. Amer. Math. Soc. 262 (1980), 429-435.



\bibitem{min}\label{MO} M. Min-Oo, E. Ruh and Ph. Tondeur, \emph{Vanishing theorems for the basic cohomology of Riemannian foliations}, J. Reine Angew. Math. 415 (1991), 167-174.


\bibitem{PJ}\label{PJ} J. S. Pak and S. D. Jung, \emph{A transversal Dirac operator and some vainshing theorems on a complete foliated Riemannian manifold}, Math. J. Toyama Univ. 16 (1993), 97-108.

\bibitem{PJ}\label{PJ1} J. S. Pak and S. D. Jung, \emph{Some vanishing theorems on complete K\"ahler foliations}, Acta Math. Hungar. 77 (1997), 15-28.

\bibitem{PR}\label{PR} E. Park and K. Richardson, \emph{The basic Laplacian of a Riemannian foliation}, Amer. J. Math. 118 (1996), 1249-1275.

\bibitem{RE}\label{RE} B. Reinhart, \emph{Differential Geometry of Foliations}, Springer-Verlag, 1983.





\bibitem{Tond1}\label{TO} Ph. Tondeur,
\emph{Geometry of foliations}, Basel: Birkh\"auser Verlag, 1997.





\bibitem{Yorozu}\label{YO1} S. Yorozu, \emph{Notes on suare-integrable cohomology spaces on certain foliated manifolds}, Trans. Amer. Math. Soc. 255 (1979), 329-341.












\end{thebibliography}
\end{document}